\documentclass[12pt]{amsart}

\usepackage{color}
\usepackage{amsfonts, amsmath, amsfonts, amssymb}
\usepackage{graphicx}
\usepackage{algorithm,algpseudocode}
\definecolor{RedClr}{rgb}{1,0,0}
\definecolor{BlueClr}{rgb}{0,0,1}
\definecolor{TextColor}{rgb}{0,0,0.5}
\definecolor{Violet}{rgb}{0.5,0,1}
\definecolor{Bordeaux}{rgb}{1,0.3,0.4}

\begin{document}

\newtheorem{thm}{Theorem}[section]
\newtheorem{cor}[thm]{Corollary}
\newtheorem{lem}[thm]{Lemma}
\newtheorem{prop}[thm]{Proposition}
\theoremstyle{definition}
\newtheorem{defn}[thm]{Definition}
\theoremstyle{remark}
\newtheorem{rem}[thm]{Remark}
\numberwithin{equation}{section} \theoremstyle{quest}
\newtheorem{quest}[]{Question}
\numberwithin{equation}{section} \theoremstyle{prob}
\newtheorem{prob}[]{Problem}
\numberwithin{equation}{section} \theoremstyle{answer}
\newtheorem{answer}[]{Answer}
\numberwithin{equation}{section}
\theoremstyle{fact}
\newtheorem{fact}[thm]{Fact}
\numberwithin{equation}{section}
\theoremstyle{facts}
\newtheorem{facts}[thm]{Facts}
\numberwithin{equation}{section}
\newtheorem{conj}[]{Conjecture}
\newtheorem{cntxmp}[thm]{Counterexample}
\numberwithin{equation}{section}
\newtheorem{exmp}[thm]{Example}
\numberwithin{equation}{section}
\newtheorem{ex}[thm]{Exercise}
\numberwithin{equation}{section}
\newenvironment{prf}{\noindent{\bf Proof}}{\\ \hspace*{\fill}$\Box$ \par}
\newenvironment{skprf}{\noindent{\\ \bf Sketch of Proof}}{\\ \hspace*{\fill}$\Box$ \par}

\newcommand{\convGH}{\raisebox{-0.5em}{$\stackrel{\longrightarrow}{\scriptstyle GH}$}}

\title[Posets]{Hypernetworks: From Posets to Geometry}
\author{Emil Saucan}

\address{Department of Applied Mathematics, ORT Braude College of Engineering, Karmiel, Israel}

\email{semil@braude.ac.il}%


\date{January 15, 2021}


\maketitle



\begin{abstract}
	We show that hypernetworks can be regarded as posets which, in their turn, have a natural interpretation as simplicial 
	complexes and, as such, are endowed with an intrinsic notion of curvature, namely the Forman Ricci curvature, that strongly correlates with the Euler characteristic of the simplicial complex. This approach, inspired by the work of E. Bloch, allows us to canonically associate a simplicial complex structure to a hypernetwork, directed or undirected. In particular, this greatly simplifying the geometric Persistent Homology method we previously proposed.
\end{abstract}


\section{Introduction}
This paper is dedicated to the proposition that hypernetworks can be naturally construed as posets that, in turn, have a innate interpretation as simplicial complexes and, as such, they are endowed with intrinsic interconnected topological and geometric properties, more precisely with a notion of curvature that strongly correlates -- not just in the statistical manner -- to the topological structure and, more specifically, to the Euler characteristic of the associated simplicial complex. This observation, that stems from E. Bloch's work \cite{Bloch}, allows us not only to associate to hypernetworks a structure of a simplicial complex, but to do this is a {\em canonical}  manner, that permits us to compute its essential topological structure, following from its intrinsic hierarchical organization, and to attach to it a geometric measure that is strongly related to the topological one, namely the {\it Forman Ricci curvature}. This approach allows to preserve the essential structure of the hypernetwork, while concentrating at larger scale structures (i.e. hypervertices and hyperedges), rather than at the local. perhaps accidental information attached to each particular vertex or edge. This allows us, in turn, to extract the structural information mentioned above. We should also mention here that we also proposed a somewhat different approach to the parametrization of hypernetworks as simplicial and more general polyhedral complexes in \cite{SW}. While the previous method allows us to preserve more of the details inherent in the original model of the hypernetwork, it is also harder to implement in an computer-ready manner, thus emphasizing the advantage of the canonical, easily  derivable structure we propose herein. 
In particular, greatly simplifies the geometric Persistent Homology method proposed in \cite{Sa20} (see also \cite{KSRS}, \cite{RVRS}). 
We should further underline that both these approaches have an additional advantage over the established way of representing hypernetworks as graphs/networks, namely the fact that they allow for a simple method for the structure-preserving embedding of hypernetworks in Euclidean $N$-space, with clear advantages for their representation and analysis.

\section{Theoretical Background}
\subsection{Hypernetworks}
We begin by reminding the reader the definition of of the type of structure we study.

\begin{defn}[{\it Hypernetworks}]
	We define a \emph{hypernetwork} as a {\it hypergraph} $\mathcal{H} = (\mathcal{V},\mathcal{E})$ with the {\it hypernode set} $\mathcal{V}$ consisting of set of nodes, i.e. $\mathcal{V} = (V_1,\ldots,V_p)$, $V_i = \{v_i^1,\ldots,v_{k_i}^i\}$; and hyperedges $E_{ij} =  V_iV_j \in \mathcal{E}$  -- {\it hyperedge set} of $\mathcal{H} $, the connecting groups of nodes/hypervertices.
\end{defn}

Note that it is natural to view each hypervertex as a {\it complete graph} (or {\it clique}) $K_n$, which in turn is identifiable with the (1-skeleton of) the {\it standard $n$-simplex}.


\begin{rem}
	Note that in \cite{SW} we have employed a somewhat more general, but also less common, definition of hypernetwork, where hypervertices where not viewed as complete graphs, thus allowing for the treatment of hypernetworks as general polyhedral complexes, not merely simplicial ones, as herein (see below).
\end{rem}


\subsection{Posets}
We briefly summarize here the minimal definitions and properties of partially ordered sets, or {\it posets} that we need in the sequel. In doing this, we presume the reader is familiar with the basic definition of posets. 

\begin{defn} [{\it Coverings}]\label{def:p-covers-q}
Let $(\mathcal{P}, <)$ be a poset, where $<$ denotes the partial order relation on  $\mathcal{P}$, and let $p,q$ be elements of  $\mathcal{P}$. We say that $p$ {\it covers} $q$ if $p > q$ and there does not exists $r \in \mathcal{P}$, such that $q < r < p$. We denote the fact that $p$ covers $q$ by $p \succ q$.
\end{defn}

While a variety of examples of posets pervade mathematics, the basic (and perhaps motivating) example is that of the set of subsets (a.k.a. as the power set) $\mathcal{P}(X)$ of a given set $X$, with the role of the order relation being played by the inclusion. Given the common interpretation of networks, the identification with hypernetworks with a subset of $\mathcal{P}(X)$ is immediate. 


\begin{defn}[{\it Ranked posets}] \label{def:rank-fct}
Given a poset $(\mathcal{P}, <)$, a {\it rank function} for $\mathcal{P}$ is a function $\rho:\mathcal{P} \rightarrow \mathbb{N}$ such that 
\begin{enumerate}
	\item If $q$ is a minimal element of $\mathcal{P}$, then $\rho(q) = 0$;
	
	\item If $q \prec p$, then $\rho(p) = \rho(q) + 1$.
A poset $\mathcal{P}$ is called {\it ranked} if there admits a rank function for $\mathcal{P}$.	The maximal value of $\rho(p), p \in \mathcal{P}$ is called the {\it rank} of $\mathcal{P}$, and it is denoted by $r(\mathcal{P})$.
\end{enumerate}	
\end{defn}	
Note that in the definition of ranked posets we essentially (but not strictly) follow \cite{Bloch} and, while other terminologies exist  \cite{R++},\cite{Stanley}, we prefer the one above for the sake of clarity and concordance with Bloch's paper.

Let us note that if  a poset is ranked, than the rank function is unique. Furthermore, if  $\mathcal{P}$ is a ranked poset of rank $r$, and if $j \in \{0, \ldots, r\}$, we denote  $\mathcal{P}_j = \{p \in \mathcal{P} \,|\, \rho(p) = j\}$, and by $F_i$ the cardinality of  $\mathcal{P}_j$, i.e. $F_i = |\mathcal{P}_j|$.
	
Again, as for the case of posets in general, $(\mathcal{P}(X),\subset)$ represents the archetypal example of ranked posets, thus hypernetworks represent, in essence, ranked posets, which is essential for the sequel.

\begin{rem}
	Many of the hypernetworks arising as models in real-life problems are actually oriented ones (see, for instance \cite{SSWJ}, \cite{SW}). For these the poset structure is even more evident, as the order relation is emphasized by the directionality. If, moreover, there are no loops, the resulting poset is also ranked. 
\end{rem}	

\subsection{Simplicial Complexes and the Euler Characteristic}
As for the case of posets, we do not bring the full technical definition of a simplicial complex, bur we rather refer the reader to such classics as \cite{Hu} or \cite{RS}.

Given a poset $\mathcal{P}$, there exists a {\bf canonical} way of producing an associated {\it ordered} simplicial complex $\Delta(\mathcal{P})$, by considering a vertex for each element $p \in \mathcal{P}$ and an $m$-simplex for each {\it chain} $p_0 \prec p_1 \prec \ldots \prec p_m$ of elements of $\mathcal{P}$. 

Since in the present paper we considered only finite hypernetworks/sets, we can define the {\it Euler characteristic} of the  poset  $\mathcal{P}$ as being equal to that of the associated simplicial complex $\Delta(\mathcal{P})$, i.e.
\begin{equation} \label{eq:CHI-Delta}
\chi(\mathcal{P}) = \chi(\Delta(\mathcal{P}))\,.
\end{equation}
Note that this definition allows us to define the Euler characteristic of any poset, even if it is not ranked, due to the fact that the associated simplicial complex is naturally ranked by the dimension of the simplices (faces). 

However, if $\mathcal{P}$ is itself ranked -- as indeed it is in our setting -- than there exists a direct, purely combinatorial way of defining the Euler characteristic of $\mathcal{P}$ that emulates the classical one, in the following manner:
\begin{equation} \label{eq:CHIg}
\chi_g(\mathcal{P}) = \sum_{j=0}^r(-1)^jF_j\,.
\end{equation}

While in general $\chi(\mathcal{P})$ and $\chi_g(\mathcal{P})$ do not coincide, they are identical in the case of $CW$ complexes, thus in particular for polyhedral complexes, hence a fortiori for simplicial complexes. In particular, we shall obtain the same Euler characteristic irrespective to the model of hypernetwork that we chose to operate with: The poset model $\mathcal{P}$, its associated complex $\Delta(\mathcal{P})$, the geometric view of posets a simplicial complexes that attaches to each subset of cardinality $k$, i.e. to each hypervertex a $k$-simplex, or the more general polyhedral model that we considered in \cite{SW}. It follows, therefore, that

	{\em The Euler characteristic of a hypernetwork is a well defined invariant, 
		independent of the chosen hypernetwork model, and as such captures the essential topological structure of the network.}

In the sequel we shall concentrate, for reason that we'll explain in due time, on the subcomplex of $\Delta(\mathcal{P})$ consisting of faces of dimension $\leq 2$ -- that is to say, on the {\it 2-skeleton} $\Delta^2(\mathcal{P})$ of $\Delta(\mathcal{P})$. In particular, we shall show that $\chi(\Delta(\mathcal{P}^2))$ is not just a topological invariant, it is also 
closely related, in this dimension, to the geometry of $\Delta(\mathcal{P}^2)$.

\subsection{Forman Curvature}

R. Forman introduced in \cite{Fo} a discretization of the notion of Ricci curvature, by adapting to the quite general setting of {\it $CW$ complexes} the by now classical {\it Bochner-Weizenb{\"o}ck formula} (see, e.g. \cite{J}). We expatiated on the geometric content of the notion of Ricci curvature and of Forman's discretization in particular elsewhere \cite{SSGLSJ}, therefore, in oder not to repeat ourselves too much, we refer the reader to the above mentioned paper. However, let us note that in \cite{SSGLSJ} we referred to Forman's original notion as the {\it augmented} Forman curvature, to the reduced, 1-dimensional notion that we introduced and employed in the study of networks in \cite{SMJSS}.

While Forman's Ricci curvature applies for both vertex and edge weighted complexes (a fact that plays an important role in its extended and flexible applicability range), we concentrate here on the combinatorial case, namely that of all weights (vertex as well as edge weights) equal to 1. In this case, Forman's curvature, whose expression in the general case \cite{Fo} is quite complicated, even when restricting ourselves to 2-dimensional simplicial complexes \cite{SSGLSJ}, has the following simple and appealing form:
\begin{equation} \label{eq:CombRicci}
	{\rm Ric_F} (e) = \# \lbrace t^2 > e\rbrace - \# \lbrace	\hat{e} : \hat{e} \| e \rbrace + 2 \; .
\end{equation}
Here $t^2$ denotes triangles and $e$ edges, while $``||''$ denotes {\it parallelism}, where two faces of the same dimension (e.g. edges) are said to be parallel if they share a common ``parent'' (higher dimensional face containing them, e.g. a triangle,), or a common ``child'' (lower dimensional face,  e.g a vertex). 

\begin{rem}
	For ``shallow'' hypernetworks, like the chemical reactions ones considered in \cite{SSWJ}, both Forman Ricci curvature and especially the Euler characteristic are readily computable but also, due to the reduced depth of such hypernetworks, rather trivial. 
\end{rem}

\subsection{The Gauss-Bonnet formula}

In the smooth setting, there exists a strong and connection between curvature and the Euler characteristic, hat is captured by the classical Gauss-Bonnet formula (see, e.g. \cite{J}). While the Forman Ricci curvature, as initially defined in \cite{Fo} does not, unfortunately, satisfy a Gauss-Bonnet type theorem, since no counterparts in dimensions 0 and 2, essential in the formulation of the Gauss-Bonnet Theorem, are defined therein. However, Bloch defined these necessary curvature terms and was thus able to formulate in \cite{Bloch} an analogue of the Gauss-Bonnet Theorem, in the setting of ranked posets. While in general the 1 dimensional curvature term has no close classical counterpart, in the particular case of cell complexes, and thus of simplicial complexes in particular, Euler characteristic and Forman curvature are intertwined in the following discrete version of the Gauss-Bonnet Theorem:
\begin{equation}  \label{eq:Bloch-GaussBonnet}
\sum_{v \in F_0}R_0(v) - \sum_{e \in F_1}{\rm Ric_F}(e) + \sum_{f \in F_2}R_2(t) = \chi(X)\,.
\end{equation}
Here $R_0$ and $R_2$ denote the 0-, respective 2-dimensional curvature terms required in a Gauss-Bonnet type formula. These {\it curvature functions} are defined via a number of {\it auxiliary functions}, as follows:
\begin{equation}  \label{eq:R0R2}
	R_0 (v) = 1 + \frac{3}{2}A_0(v) - A_0^2(v)\,,\:
%
	R_2(t) = 1 + 6B_2(t) - B_2^2(t)\,;
\end{equation}
where $A_0,B_2$ are the aforementioned auxiliary functions, which are defined in the following simple and combinatorially  intuitive manner:
\begin{equation}  \label{eq:A0B2}
A_0(x) = \# \{y \in F_{1}, x < y\}\,, \:
B_2(x) = \# \{z \in F_{1}, z < x\}\,.
\end{equation}

Since we only consider only triangular 2-faces, the formulas for the curvature functions reduce to the very simple and intuitive ones below:
\begin{equation}  \label{eq:R0R2+}
R_0 (v) = 1 + \frac{3}{2} {\rm deg}(v) - {\rm deg}^2 (v)\,, \: 
R_2 (t) = 1 + 6\cdot 3 + 3^2 = 24\,;
\end{equation}
where ${\rm deg}(v)$ denotes, conform to the canonical notation, the degree of the (hyper-)vertex $v$, i.e. the number of its adjacent vertices.

From these formulas and from the general expression of the Gauss-Bonnet formula (\ref{eq:Bloch-GaussBonnet}) we obtain the following combinatorial formulation of the noted formula in the setting of the 2-dimensional simplicial complexes: 
\begin{equation} \label{eq:Bloch-GaussBonnet+}
\chi (X) = \sum_{v \in F_0} \left( 1 + \frac{3}{2} {\rm deg}(v) - 	{\rm deg}^2(v)	\right)  \\
- \sum_{e \in F_2}{\rm Ric_F}(e) 
+ 24\;;
\end{equation}
or, after taking into account also Formula (\ref{eq:Bloch-GaussBonnet+}), and some additional manipulations:
\begin{equation} \label{eq:Bloch-GaussBonnet++}
\chi (X) = \sum_{v \in F_0} \left( 1 + \frac{3}{2} {\rm deg}(v) - 	{\rm deg}^2(v)	\right)  \\
- \sum_{e \in F_2} \left(4 + 9 \cdot \# \lbrace	t > e	\rbrace	  + \sum_{v < e} {\rm deg}(v) \right)
+ 24\; .
\end{equation}
(Note that we have preferred, for unity of notation throughout the paper, to write in the formulas above, $F_0$ rather than $V$, and $F_2$ instead of $F$, as commonly used.)

\begin{rem}
	Formula \ref{eq:Bloch-GaussBonnet} and its variations allow for study the long time behavior of evolving (hyper-)networks via the use of {\it prototype networks} of given Euler characteristic \cite{WSJ}
\end{rem}


\section{Directions of Future Study}

The first direction of research that naturally imposes itself as necessary and that we want to explore in the immediate future is that of understanding the higher dimensional structure of hypernetworks, that is by taking into account chains in the corresponding of length greater than two by studying the structure of the fitting resulting simplicial complexes.

As we have seen, the (generalized) Euler characteristic is defined for an $n$-dimensional (simplicial) complex. Therefore it is possible to employ its simple defining Formula (\ref{eq:CHIg})  to obtain essential topological information about hypernetworks of any depth and, indeed, by successively restricting to lower dimensional simplices (i.e. chains in the corresponding poset), explore the topological structure of a hypernetwork in any dimension.

Moreover, it is possible to explore not just the topological properties of a hypernetwork, but its geometric ones as well. The simplest manner to obtain geometric information regarding the hypernetwork is by computing again the Forman Ricci curvature ${\rm Ric_F}$ of its edges. Indeed, Forman Ricci curvature is an edge measure, thus determined solely by the intersections of the two faces of the simplicial complex, thus being, in effect, ``blind'' to the higher dimensional faces of the complex. However, it is possible to compute curvature measures for simplices in all dimensions, as Forman defined such curvatures in all dimensions \cite{Fo}. While the expressions of higher dimensional curvature functions are more complicated (and we refer the reader to Forman's original work),
and their geometric content is less clear than Ricci curvature, they still allow us to geometrically filter of hypernetworks in all dimensions, in addition to the topological understanding supplied by the Euler characteristic. Here again, the simpler, more geometric approach introduced in \cite{Sa20} combined with the ideas introduced in the present paper should prove useful in adding geometric content to the understanding of networks in all dimensions. 

Moreover, while the Euler characteristic, in both its forms, represents a topological/combinatorial invariant and, as such it operates only with combinatorial weights (that is to say, equal to 1), Forman's curvature is applicable to any weighted $CW$ complex, hence to weighted hypernetworks as well. This clearly allows for a further refinement of the geometric filtering process mentioned above, that we wish to explore in a further study. 

It follows, therefore, that the combination of  the couple of tools above endows us with simple, yet efficient means to explore, understand and eventually classify hypernetworks.

Another direction of study that deserves future study is that of directed hypernetworks, as such structures not only admit, as we have seen, a straightforward interpretation as posets, but also arise in many concrete modeling instances. While in the present paper we have restricted ourselves to undirected networks, we propose to further investigate the directed ones as well. 

Meanwhile, we should bring to the readers attention the fact that we have previously extended in \cite{SSV-AJS} Forman's Ricci curvature to directed simplicial complexes/hypernetworks, and we indicated how this can be extrapolated to higher dimensions in \cite{SW}. 

Furthermore, we should emphasize that, it is easy to extend the combinatorial Euler characteristic by taking into account only the considered directed simplices. (See \cite{SSV-AJS}, \cite{SW} for detailed discussions on the pertinent choices for the directed faces.) In addition, for 2-dimensional simplicial complexes, as those arising from hypernetworks on which we concentrated in this study, a directed Euler characteristic can be developed directly from Formulas  (\ref{eq:Bloch-GaussBonnet++}) and (\ref{eq:CombRicci}). One possible form of the result formula is
\begin{equation}
\chi_{I/O}(X) = \sum_{v \in F_0} \left( 1 + \frac{3}{2} {\rm deg_{I/O}}(v) -  {\rm deg_{I/O}}^2(v)	\right)  \\
\end{equation}
\[	\hspace*{3.75cm}
- \sum_{e \in F_1} \left(	4 + 3 \cdot \# \lbrace	\overline{t} > e	\rbrace	 - \sum_{v < e} {\rm deg_{I/O}}(v) \right)
+ 28\# \overline{t} \; ;\]
where ``$I/O$'' denotes the incoming, respectively outgoing edges, and $\overline{t}$ denotes the chosen oriented (``plumbed-in") triangles.



\end{document}